\input amssym.def
\input amssym
\def\idl{$\frak I$\ }
\def\gmw{$\cal G_{\omega} $ }
\def\gmu{$\cal G_{\cal U}$ }
\def\for{$\Bbb P (\cal U )$ }
\documentstyle{article}
\begin{document}
\title{Combinatorics on Ideals and Axiom A}
\author{James D. Sharp}
\footnotetext{This paper forms part of
the author's dissertation, written under the direction of Simon
Thomas}
\date{1/28/1993}

\maketitle

\section{Introduction}
The following notion of forcing was introduced by Grigorieff \cite{G}:
Let $\frak I \subset \omega$ be an ideal, then $\Bbb P$ is the set of all functions
$p:\omega \rightarrow 2$ such that $dom(p)\in \frak I$. The usual
Cohen forcing corresponds to the case when $\frak I$ is the ideal of
finite subsets of $\omega$. In \cite{G} Grigorieff proves that if \idl is the dual of a
p-point ultrafilter, then $\omega_{1}$ is preserved in the generic
extension. Later, when Shelah introduced the notion of proper forcing,
many people observed that Grigorieff forcing was proper. One way of
proving this is to show that player II has a winning
strategy in the game $\cal G_{\omega}$ for $\Bbb P$ (see \cite{J},
page 91.)

The notion of Axiom A forcing was introduced by Baumgartner \cite{BM}.
If $\Bbb P$ satisfies Axiom A, then player II has a winning strategy
in the game \gmw and thus is proper. Indeed, most of the naturally
occurring proper notions of forcing are Axiom A (e.g. Mathias or Laver
forcing). Thus it is natural to ask whether or
not Grigorieff forcing satisfies Axiom A. The main result of this
paper is a negative answer to this question. We will prove this by
introducing another game \gmu and showing that if $\Bbb P$ were Axiom A
then player II would have a winning strategy in this game. We will then prove
that the game \gmu is undetermined.

\section{Definitions and Preliminaries}
Throughout this paper $\cal U$ will denote a p-point ultrafilter and
\idl will denote the dual ideal. We let ${[X]}^{< \omega}=\{s\subseteq
X \mid |s|<\omega\}.$ And we will denote ordinal names by
$\dot{\alpha}$. We let $Seq(X)$ denote the set of finite sequences
of elements of X. If $s=\langle x_{0}, \ldots , x_{n}\rangle\in Seq(X)$
and $y\in X,$ then $s*\langle y\rangle = \langle x_{0},\ldots , x_{n},
y\rangle\in Seq(X).$
\bigskip\pagebreak
\\DEFINITION. The game \gmu is for two players playing alternatively.
\\Player I plays a partition of $\omega, \{I_{n}\mid n\in\omega\},$ such that
for all n, 
$I_{n}\in\frak I$ and player II plays finite subsets $F_{n}\subset I_{n}$.
Player II wins iff 
\begin{equation}
\bigcup_{n\in \omega}F_{n} \in \cal U.
\end{equation}

Note that this definition makes sense for an arbitrary non-principal
ultrafilter $\cal U$. However, if $\cal U$ is not a p-point, then
player I clearly has a winning strategy. \\
DEFINITION. A partial order $\Bbb P$ is said to satisfy Axiom A if there is a
collection of partial orders $\{\leq_{n} \mid n\in\omega\}$ of $\Bbb P$
satisfying:
\newcounter{ax.a}
\begin{list}{\roman{ax.a}}{\usecounter{ax.a}}
\item[i)]  $p\leq_{0} q$ implies $ p\leq q$
\item[ii)] $p\leq_{n+1} q$ implies $p\leq_{n} q$
\item[iii)]  if $ \langle p_{n}\mid n\in\omega\rangle$ is a sequence such that
\\$p_{0}\geq_{0} p_{1}\geq_{1}
p_{2}\geq_{2}\ldots\geq_{n-1}p_{n}\geq_{n}\ldots$ (called a fusion
sequence),\\ then there is a $q\in\Bbb P$ such that for all n,
$q\leq_{n}p$
\item[iv)]  for all $p\in\Bbb P$, for all $n\in\omega$, and for all ordinal
names $\dot{\alpha}$, there exists $q\leq_{n} p$ and a countable set B
such that $q\Vdash \dot{\alpha}\in B.$ 
\end{list}
\bigskip
DEFINITION. Grigorieff Forcing
\\${\Bbb P} ({\cal U}) =\{ p:\omega \rightarrow 2\mid dom(p)\in
\frak I \}$ where $q\leq p$ iff $q\supseteq p$.
\medskip
\\The main results in this paper are the following two theorems.
\newtheorem{theorem}{Theorem}
\begin{theorem}
Gregorieff forcing does not satisfy Axiom A.\label{th1}
\end{theorem}

\begin{theorem}
The game \gmu is not determined.\label{th2}
\end{theorem}

\section{Grigorieff Forcing and \gmu}

In this section we will prove the following Lemma:
\newtheorem{lemma}{Lemma}
\begin{lemma}
If $\Bbb P (\cal U )$ satisfies Axiom A, then player II has a winning
strategy in the game \gmu .\label{axa}
\end{lemma}
PROOF. Suppose  \for \/satisfies Axiom A.

\newtheorem{claim}{Claim}
\begin{claim}
Let $p\in\Bbb P (\cal U)$, $n\in\omega$, and $I\in$\idl, then there exists
$q\in$\for\/ such that 
\begin{equation}q\leq_{n} p\; and\; \mid I\setminus
dom(q)\mid < \omega.\end{equation}
\end{claim}
PROOF OF CLAIM 1: 
Let p, n, and I be given.
We may assume I is infinite.
Let \{$Y_{\alpha}\mid \alpha < 2^{\omega}$\} be an enumeration of $\wp
(I)$ and set $\dot{\alpha}=\{\langle \alpha ,
p_{\alpha}\rangle \mid \alpha < 2^{\omega}\}$ where
$p_{\alpha} : I \rightarrow 2$ is the characteristic function of
$Y_{\alpha}$. Since \for satisfies Axiom A, by (iv) there's
$q\in$ \for and a countable set B, such that $q\leq_{n} p$   and $q\Vdash
\dot{\alpha}\in B.$ But then $\mid I\smallsetminus dom(q)\mid <
\omega$ as required. $\square$

Now we describe a winning strategy for player II in the game \gmu.
Suppose player I plays $I_{0}$ at the $o^{th}$ move, then player II sets
$p_{0}=\chi_{I_{0}}$ and plays $F_{0}=\emptyset.$ After the $n-1^{st}$
turn, player I has played $I_{0},\ldots , I_{n-1}$ and player II has
played $F_{0}, \ldots , F_{n-1}$ and chosen $p_{0}{\geq}_{0}
p_{1}{\geq}_{1}\ldots {\geq}_{n-2} p_{n-1}.$ At the $n^{th}$ move player
I plays $I_{n}.$ Then by the claim there exists $p_{n}\in$ \for such
that $p_{n}{\leq}_{n-1} p_{n-1}$ and $\mid I_{n}\setminus
dom(p_{n})\mid < \omega$. Thus player II can play
$F_{n}=I_{n}\setminus dom(p_{n}).$ At the end of the game $\langle p_{n}\mid
n\in\omega\rangle$ forms a fusion sequence, and it follows from Axiom
A (iii), that $\cup
dom(p_{n})\in \frak I$. And thus player II wins.
\rightline{$\blacksquare$}
\\Notice that lemma~\ref{axa} and theorem~\ref{th2} imply theorem~\ref{th1}.

\section{The Game $G_{\cal U}$ is undetermined}
We shall prove the theorem as two lemmas.
\begin{lemma}
Player I does not have a winning strategy in the game \gmu.
\end{lemma}
PROOF. Suppose, by way of contradiction, that player I has a winning
strategy. Let $\sigma :Seq({[\omega]}^{<\omega})\rightarrow\frak I$ be
player I's winning strategy. Then ran($\sigma$) is countable. Let
$\langle Z_{n}\mid n \in \omega\rangle$ be an enumeration of
ran($\sigma$) such that $\sigma (\emptyset )=Z_{0},$ let
$J_{0}=Z_{0}$ and for $n\geq 1$ let $J_{n}=Z_{n}\setminus\cup_{j <
n}J_{j}$. Since $\cup ran(\sigma)=\omega$, the set $\langle J_{n}\mid
n\in\omega\rangle$ forms a partition of $\omega$. Thus there exists
$Y\in{\cal U}$ such that $\mid Y\cap J_{n}\mid < \omega$ and hence
$\mid Z_{n}\cap Y\mid < \omega$ for all
$n\in\omega.$
Now consider the following game where player I plays by $\sigma$ and
player II plays 
\begin{equation}
F_n=Y\cap\sigma(\emptyset ,F_0,\ldots ,F_{n-1})
\end{equation}
Notice that for all n, there exists some m such that
\begin{equation}
F_{n}=Y\cap
Z_{m}. \end{equation}
Clearly $\cup F_{n}=Y\in {\cal U}$, so player II wins. $\blacksquare$
\medskip
To prove lemma \ref{plII} we shall need the following definitions and
results from Grigorieff \cite{G}.\medskip
\\DEFINITION.
\newcounter{tree}
\begin{list}{\roman{tree}}{\usecounter{tree}}
\item[i)] $A\subset Seq([\omega]^{< \omega})$ is a p-tree if it is
non-empty and closed under taking initial segments.
\item[ii)] If $s\in A$, the ramification of A at s, denoted $R_{A}(s)$, is
the set of all $a\in {[\omega]}^{< \omega}$ such that $s\ast\langle
a\rangle\in A$.
\item[iii)] A is an ${\cal I}$-p-tree if for any $s\in A$, there exists $X\in
\cal U$ such that $[X]^{< \omega}\subseteq R_{A}(s)$.
\item[iv)] H is an $\cal I$-p-branch of A if it is a branch such that $\cup
\{H(n)\mid n\in\omega\}\in\cal U$.
\end{list}

\begin{theorem}
{\rm (Grigorieff \cite{G})} Every $\cal I$-p-tree has an $\cal I$-p-branch.
\end{theorem}
\begin{lemma}
Player II does not have a winning strategy in the game $G_{\cal U}$.
\label{plII}\end{lemma}
PROOF. Suppose, by way of contradiction, that player II has a winning
strategy $\mu$.
Player I will construct a ``tree of games'' such that along each
branch player II plays by $\mu$,  but there will be a branch such
that player I wins the corresponding game. First we need the following
technical result.
\begin{claim}
For all $n\in\omega$ and for any sequence $I_{0}, I_{1},\ldots , I_{n}$
of pairwise disjoint elements of \idl, there exists $Y\in \cal U$,
$Y\subseteq (\omega\setminus \cup\{I_{n}\mid i\leq n\})$
satisfying:
\begin{equation}
(\forall t\in [Y]^{<\omega})(\exists I\in {\frak I})[(t\subset
I\subset(\omega\setminus \bigcup_{i\leq n}I_{i}))\wedge (\mu
(I_{0},\ldots I_{n}, I)\cap t = \emptyset)].\label{eq} \end{equation}

\end{claim}
PROOF OF CLAIM. Let n, $I_{0},\ldots , I_{n}$ be given. Suppose the
claim is false. Then we can define, by induction on $m\in\omega$, a
sequence of pairs $\langle X_{m}, s_{m}\rangle$ such that
\newcounter{xm.sm}
\begin{list}{\roman{xm.sm}}{\usecounter{xm.sm}}
\item[i)] $X_{0}=\omega\smallsetminus\bigcup_{i\leq n}I_{i}$
\item[ii)] for all $m\geq 1,\; X_{m}=X_{m-1}\setminus s_{m-1}$
\item[iii)] $(\forall I\in{\frak I})(s_{m}\subset I\subset X_{0})\rightarrow
\mu (I_{0},\ldots ,I_{n}, I)\cap s_{m}\neq\emptyset$
\end{list}
If we set $J_{0}=\bigcup_{i\in\omega}s_{2i}$ and
$J_{1}=\bigcup_{i\in\omega}s_{2i+1}$ then one of these must be in \idl
, say $J_{0}\in\frak I$. But then for all $i\leq\omega ,\; \mu
(I_{0},\ldots , I_{n}, J_{0})\cap s_{2i}\neq\emptyset$, so $\mid
\mu (I_{0},\ldots , I_{n}, J_{0})\mid = \omega$, contradicting the
assuption that $\mu$ is a winning strategy for player II.$\square$

Player I constructs a tree $A\subset Seq(\{\langle I, \;F, s\rangle\mid
I\in{\frak I}, F,s\in [I]^{< \omega}\})$ by induction on the height
n of a node such that, if $\langle\langle I_{0}, F_{0}, s_{0}\rangle
,\ldots ,\langle I_{n}, F_{n}, s_{n}\rangle\rangle\in A$ then
\newcounter{gtree}
\begin{list}{\roman{gtree}}{\usecounter{gtree}}
\item[i)] $I_{i}\cap I_{j} = \emptyset\; if\; i\neq j$
\item[ii)] $(1\leq i\leq n)(F_{i}=\mu (I_{1},\ldots , I_{i}))$
\item[iii)] $(0\leq i\leq n)(s_{i}\cap F_{i} = \emptyset)$
\item[iv)] $I_{0}=F_{0}=s_{0}=\emptyset .$
\item[v)] The projection map $\pi :\langle\langle I_0,F_0,s_0\rangle
,\ldots ,\langle I_n, F_n, s_n\rangle\rangle\rightarrow\langle
s_0,\ldots ,s_n\rangle$ is an injection.
\end{list}

Suppose we've decided which $B=\langle \langle I_{0}, F_{0},
s_{0}\rangle,\ldots , \langle I_{n}, F_{n}, s_{n}\rangle\rangle\in A.$
Fix such a B. By the claim, there exists $Y\in {\cal U}$ satisfying
property~\ref{eq}. For each $t \in {[Y]}^{< \omega}$, choose $I_{t}$
satisfying property~\ref{eq}. Let $R_{A}(B)=\{\langle I, F, s\rangle\mid s\in {[Y]}^{<
\omega}, I=I_{s}, F=\mu (I_{0},\ldots , I_{n}, I_{s})\}.$

Let $\tau\subset Seq({[\omega]}^{<\omega})$ be the tree obtained from
A via the projection map $\pi$ 
Then, by the construction, $\tau$ is an $\frak I$-p-tree, and hence has an
$\frak I$-p-branch H. 
Let $\langle\langle I_n,F_n,H(n)\rangle\mid n<\omega\rangle$ be the
corresponding branch in A. 
Then, since $\cup \{H(n)\mid n\in \omega \}\in \cal U$, this is a play
of $G_{\cal U}$ which player I wins. 
But player II has used the strategy $\mu$ in this play, which
contradicts the assuption that $\mu$ is a winnig
strategy for player II.$\blacksquare$

\noindent{\bf CURRENT ADDRESS:}\\
James D. Sharp\\
Department of Mathematics\\
Rutgers University\\
New Brunswick, NJ 08903\\
email: jsharp@math.rutgers.edu
\end{document}